\documentclass[11pt]{amsart}
\usepackage{graphicx}
\usepackage{amssymb}
\usepackage{epstopdf}
\usepackage{amsmath,amscd}
\usepackage{amsthm}
\usepackage{verbatim}
\usepackage{stmaryrd}
\usepackage{caption,subcaption}
\usepackage{sidecap}
\usepackage[dvipsnames]{xcolor}

\usepackage[hyphens]{url}
\RequirePackage[colorlinks,citecolor=blue,urlcolor=blue]{hyperref}

\theoremstyle{plain}
\newcommand\CC{{\mathbb C}}

\newcommand\HH{{\mathbb H}}
\newcommand\RR{{\mathbb R}}
\newcommand\ZZ{{\mathbb Z}}

\newcommand\PP{{\mathbb P}}

\newcommand\sG{{\mathcal G}}

\newcommand\eps{\epsilon}
\newcommand\what{\widehat}

\newcommand\q{\quad}

\newcommand\oo{\infty}

\newcommand\pc{p_{\text{\rm c}}}
\newcommand\sL{\mathcal{L}}

\newcounter{mycount1}\newcounter{mycount2}\newcounter{mycount3}\newcounter{mycount}

\numberwithin{equation}{section}
\numberwithin{theorem}{section}
\numberwithin{figure}{section}

\title[Dominic Welsh (1938--2023)]{Dominic Welsh\\ (1938--2023)}
\author{Geoffrey R. Grimmett}
\address{Statistical Laboratory, Centre for
Mathematical Sciences, Cambridge University, Wilberforce Road,
Cambridge CB3 0WB, UK} 
\email{grg@statslab.cam.ac.uk}
\urladdr{\url{http://www.statslab.cam.ac.uk/~grg/}}

\begin{document}
\begin{abstract}
This biographical and scientific memoir of Dominic Welsh includes summaries of his important contributions 
to probability and combinatorics.
With John Hammersley, he introduced first-passage percolation, and in so doing 
they formulated and proved the first 
subadditive ergodic theorem. 
Welsh has numerous results in matroid theory, and wrote the first monograph on the topic. 
He worked on computational complexity and particularly the complexity of
computing the Tutte polynomial.  He was an inspirational teacher and advisor who
helped to develop a community of scholars in combinatorics.
\end{abstract}

\date{\today}

\maketitle

\centerline{\includegraphics[width=0.5\hsize]{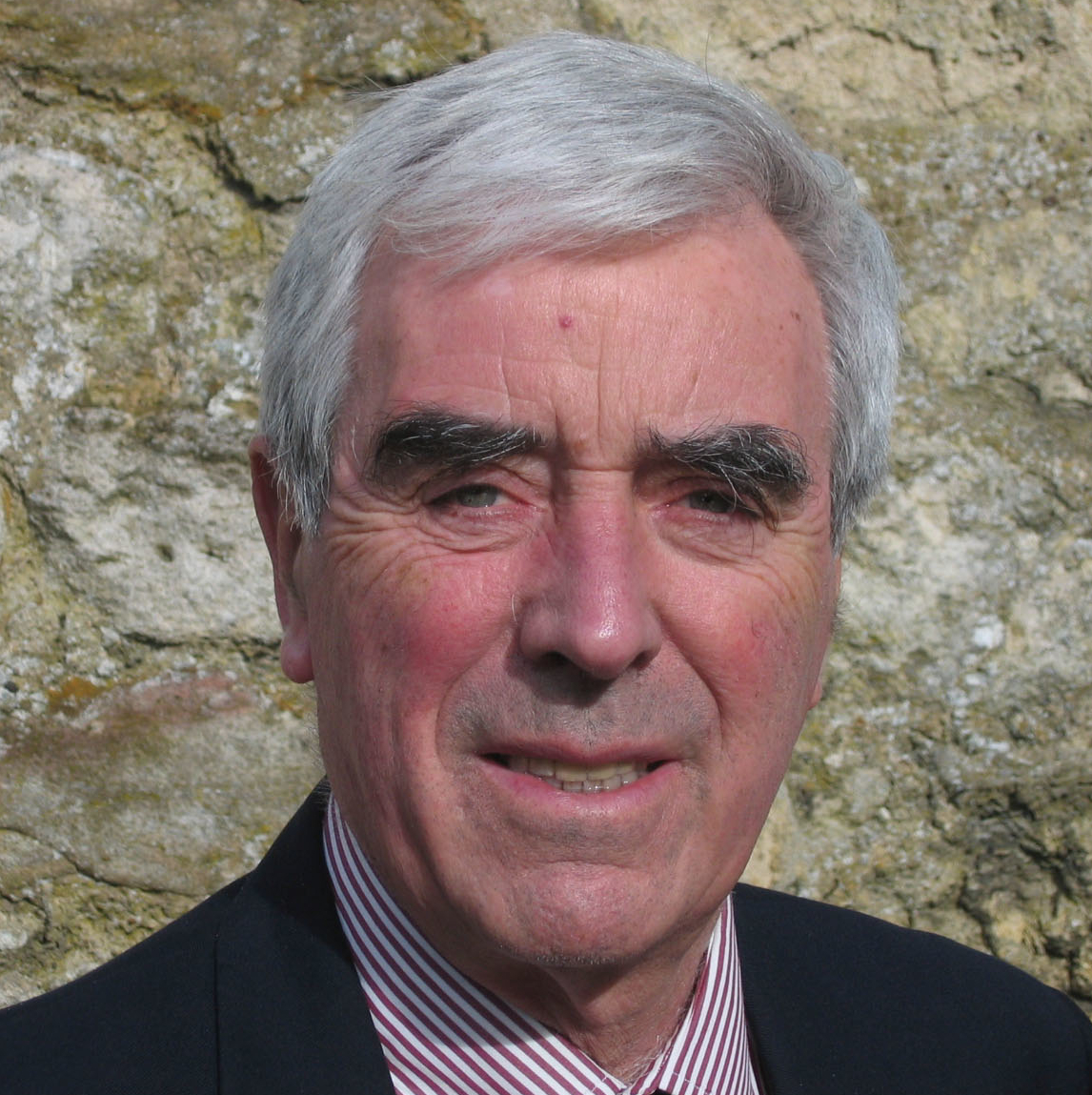}}
\bigskip

If the principal qualities of an academic are associated with  teaching, research, and collegiality, Dominic Welsh
was an outstanding example to all. He combined excellence as tutor 
and supervisor over nearly 40 years with a distinguished
research record in probability and discrete mathematics, while creating
a community of younger people who remember him with love and respect.  

\section{Early life}

James Anthony Dominic Welsh was born on 29 August 1938 in Port Talbot, the eldest of four children of parents
Teresa O'Callaghan and James Welsh. They addressed him as Dominic, and 
it was thus that he became known to all.

Dominic was born into an extended
family of Irish Catholic origins living in and around Swansea. 
Numerous members of the family including both his parents
were teachers,
and his father was Headmaster of St David's School, Swansea\footnote{St David's Roman Catholic School
is said to have been built to “save the Irish and Belgian children of the district from the risk of being 
brought up in ignorance of their religion, if not losing it altogether”, \cite{std}.}.
Education was naturally considered a priority for their four children,
though teaching was not the only profession of the family; 
one of Dominic's grandfathers worked at the Port Talbot steelworks,
and his uncle Jim (Bennett), a Dubliner by origin, was the brewer at the historic Buckley's Brewery
in Llanelli. 
All four of Teresa and James's children followed their parents' example in one way or another.
Their second son, David (b.\ 1939), studied in London and Oxford, 
and became Professor of Classics at the University of Ottawa;
Mary (b.\ 1942)  studied in Swansea, became a teacher, and moved to Japan before starting 
a successful EFL (English as a Foreign Language) business; 
Teresa (b.\ 1952) graduated with a first-class degree in mathematics 
at Swansea, and remained in Port Talbot as a teacher and later a school governor.

Dominic's earliest memories
were of the convivial home of his (maternal) Gran in Castle Street,
Port Talbot, where his aunt May married Jim Bennett to an accompaniment of 
Irish dancing, music, and songs. With the start of war, his father's school was evacuated from Swansea, and the family
moved with the newly born Dave to Ammanford. Later he returned to Port Talbot,
where he lived until a bombing raid struck a nearby house.

As a child Dominic developed a talent for rugby; a love of
sports, games, and other competitive activities stayed with him throughout his life. 
He played for the Glamorgan and District Schoolboys Team  against London at Twickenham in 1957,
capitalising on this highlight by scoring a try between the famous posts. 
He was invited to trial for the Wales Under-18 XV, but
was disappointed to be disqualified as being three days too old for inclusion in the team.
And in any case he was discouraged by the line of dentures hanging on hooks in the changing room.

For his secondary education, Dominic attended Bishop Gore Grammar School, where in due
course he became Head Boy. 
His mathematical
abilities were evidently exceptional, and he was the first boy ever to move from the school to Oxford University.
Indeed he was offered places at both Oxford and Cambridge, but the latter would have required him to first spend 
two years in National Service.  

\section{Oxford and Merton}

Dominic went up in 1957 to read Mathematics at Merton College, Oxford, where he was awarded
a scholarship (an Exhibition). 
The College became his intellectual home for evermore.
Following his BA in 1960, he visited North America for the first time with a Fulbright award to Carnegie Mellon University,
and thence to Bell Telephone Laboratories in Murray Hill. He returned to Oxford in 1961 as a postgraduate,
and was awarded a NATO studentship in 1962.
He took his D.Phil.\ in 1964 under John Hammersley, was appointed Junior Lecturer in the Mathematical Institute in 1963,
and was elected to a Tutorship and Fellowship at Merton in 1966.  
Within the University, he was promoted to a Readership in 1990 and a personal Chair in 1992,
while retaining his Fellowship at Merton.
He attained the retirement age of 67 in 2005 becoming Emeritus Fellow and Professor.
He died in Oxford on 30 November 2023 after a period of illness. 
His funeral and Requiem Mass
took place in the College Chapel on 16 December 2023, followed on 1 June 2024 by a memorial service. 

His 39 years as a Fellow  
have been  exceeded by only few in recent times, including by his colleague Philip Watson (Fellow, 1950--1993),
with whom Dominic shared the privilege of teaching and advising the many 
mathematics undergraduates at the College during their years in common. No one has contributed more
than Dominic to the marked rise in reputation of the College in mathematics over the last 60 years.

Quite apart from his family life (see Section \ref{sec:fam}), 
Dominic was a very busy person indeed. As a so-called CUF lecturer,
he had a tutorial load of 12 hours per week, plus a lesser lecturing load.
Administration had a lighter touch in the 1960s and \rq 70s than now, but Dominic was never a shirker and he undertook
his share of duties, including the College-appointed posts of Principal of the Postmasters [PoP]
(1970--73), Pro-Proctor (1979--80), Sub-Warden (1982--84), and Wine Steward (1998--2002).

Within the University,  he served as Chair
of the Faculty (1976--78), of the Board of Mathematical Sciences (1984--86), and of the 
Mathematical Institute (1996--2001). 
His term as Institute Chair was eminently successful, though he
was sometimes uncomfortable with the degree of influence and control available to him in the role. 

Dominic was elected to the London Mathematical Society in 1969, and served as 
Publications Secretary and Member of Council for 1972--1976. 
He was a member of the editorial boards of numerous leading journals, 
and was a popular lecturer at international conferences and workshops.

Dominic's strong sense of responsibility was invariably leavened with humour and good sense,
and his deepest loves were people and research.  
He was very good with young people, and he developed a
rare and enviable affinity with many of his undergraduate and postgraduate students. 

Dominic was an academic of a now rare breed. He was devoted and generous to his students, 
and he welcomed them into his home and family
for nourishment, conversation, and other (usually competitive)
activities. He taught up to 15 hours per week during term, and covered a wide range of
subjects including the entire first- and second-year syllabuses. He was active over more than 50 years
at a high level in mathematical research; he supervised 28 D.Phil.\ students,
and he wrote numerous research papers and books with more than 50 collaborators. 
He undertook a full load of administration in the College and University.
All this he achieved with charm and diligence. 

Following his retirement from Oxford, Dominic and Bridget went on tour, 
living (recalls Bridget) out of seven suitcases for two years.
They travelled first to New Zealand for a three-month visit to the University of Canterbury, and then
to Bath for six weeks at the beginning of an association with the Heilbronn Institute in Bristol 
that was to last until 2012.  
They particularly enjoyed their year in Barcelona where Dominic was visiting the Technical University.
On 17 August 2007 they had their first experience of owner-occupancy on moving into
their new home in Charlbury Road. 

\section{Family life}\label{sec:fam}

Dominic and Bridget met on Dominic's 24th birthday in 1962. It was a \emph{coup de foudre}.
They married on 10 July 1965 and made their first home together in a Woodstock cottage
backing onto the Blenheim Estate. Their three sons 
arrived in arithmetic progression: James (b.\ 1967),
Simon (b.\ 1969), John (b.\ 1971). The family came to know many of Dominic's students,
initially through childsitting, and progressing to other activities. 

\emph{Dominic's family has 
written} (in \cite{fam}):
\lq Many of his colleagues and students went on to be lifelong friends. Dominic and Bridget would host 
long laughter-filled Sunday lunches, first at their college home in Kybald Street then later in Rose Lane, 
where future mathematical geniuses would sit side by side with a trio of grubby-faced mischievous kids.'
[Some readers may cavil at the word \lq future'.]

Two aspects of Dominic's personality merit mention. He was born into a Catholic family 
and was committed to Roman Catholicism throughout his life. He was 
guided by his religion  in moral matters, but without judgement of  others.  
He attended the Catholic Chaplaincy regularly 
before moving to the Oxford Oratory.

Secondly, he was an avid and competitive sports and games player. Dinner would often be followed by 
a game, usually designed with children in mind but played with the seriousness of a rugby final at Cardiff Arms Park. 

\emph{David Stirzaker has written}:
\lq Dominic's  pastoral care of students extended to playing both tennis and
real tennis with them on the Merton courts.  He was fond of describing real
tennis as \lq chess on the run', $\dots$ or was it \lq chess on wheels',   I
forget.   He enjoyed beating me at lawn tennis.
He used to invite students round to his house to play group board
games.  I remember enjoying some very competitive sessions of Diplomacy
in the 1960s, which revealed a rich vein of deviousness and
cunning in his character. 
 Students and colleagues were invited to watch rugby with him.  
 I well recall being at his house in Kybald St to see the Barbarians
beat New Zealand in 1973.  His joy at Gareth Edwards' try was very 
infectious.\footnote{David Williams is quoted in \cite{DW} as saying
  ``[Maths is] like watching the famous match between the Barbarians and the All Blacks. Everyone who saw that game was carried away with the excitement. But it does not compare with the excitement of mathematics.''}\rq\
[His family considered him a good loser but a horrible winner.]

Disaster struck in 1990 when John drowned in an accident in Brisbane.
To quote again from the Welsh family:
\lq They were happy times, but no life is without pain. And in 1990, Dominic and Bridget suffered the worst pain imaginable, when their beloved eighteen-year-old son, John, died while travelling in Australia. No one ever fully recovers from the loss of a child, but eventually Dominic and Bridget somehow found the strength to keep on living, loving and laughing.\rq

In later life, Bridget and Dominic travelled the world together and enjoyed walking holidays.

\section{Tributes, reminiscences}

\emph{Adrian Bondy has written}:
\lq As my doctoral supervisor, Dominic proposed to me a good number of thesis topics. I worked on matroids for a while, coming up with a variety of matroids based on graphs (which all turned out to be derivable from graphic matroids by simple operations). Finally, I latched onto the reconstruction problem, which remains unsettled. 

Dominic was always bubbling over with ideas. His enthusiasm was contagious. Just a few years my senior, he was friendly and approachable, and I came to know him and his family quite well. It was an ideal environment for a research student.\rq

\emph{James Oxley has written}: 
\lq Dominic modelled for me not only how to do research but also how to
bring out the best in students. The first time I met him was in the Porters' Lodge
at Merton College just after I had arrived in Oxford.
Seeing him there, I tentatively approached him with \lq\lq Dr Welsh?"
\lq\lq Dominic" was his immediate response, quickly followed by an invitation
to afternoon tea with his family. His friendly, welcoming, and informal manner
made interacting with him great fun.

Dominic's generosity of spirit shone through in our exchanges. After
each of us had published books on matroid theory,
he apologised to me because he was worried that the 2010 Dover reprint
of his 1976 book may hurt the sales of my book.
Both Dominic's sage guidance and our deep and enduring friendship enriched
my life immensely.'

\emph{Andrew Wiles has written}:
\lq I first met Dominic at my undergraduate interview in Merton. I don't remember much of it except for his genuine
warmth. He asked me at one point to explain epsilon/delta arguments and I thought of making a stab at it, but 
decided it was better not to make a hash of it. In truth I had really never learned it. So after a thirty second silence 
I just said that I don't know and he brushed over it and put me at my ease. It was the same when I was his student 
in the first year. He accepted that I had strong tastes in mathematics and did not complain that my performance was 
much better in some areas than in others. 

For the rest of his life he and Bridget would always welcome me back to 
Oxford in the summers and especially in the year 1988--90. When I finally moved back for good in 2011 I moved 
into what had been his office and the memories of earlier visits came back with it. Not everything in the all male 
Merton world of my undergraduate days would I wish to remember but Dominic's presence was certainly the best 
part of it and I will always cherish his friendship and kindness then and afterwards.'

\emph{Neil Loden has written}:
\lq I was sad to learn of the death of Dominic Welsh whom I remember fondly in his role as the PoP.  
After a particularly debauched dinner of the 1311 Club in 1969, 
Alan [Harland], Chris [Hewitt] and I  blocked up the door into his rooms by building a wall of bricks (taken 
from the works going on in Front Quad) from floor to ceiling, which kept him out of his room for the whole 
of the next morning, while we were nowhere to be seen, being still in bed nursing hangovers. We thought at 
the time that it was very funny, and even he didn't seem to mind unduly.\rq
 
 \emph{David Stirzaker has written}:
\lq Dominic's pupils generally found him to be one of
the most engaging and entertaining of all their lecturers.  More than other lecturers, 
he conveyed the impression that he was enjoying it
and found it fun to be talking to us.    I recall him remarking of some
theorem that, as it stood, it was so hedged about with technical conditions
and restrictions as to be about as useful as a barren pear tree.   
This was quite refreshing after the arid dustiness of many Oxford lectures.
He was even more outstanding as a tutor, because he had the twin gifts of
seeing precisely where the student's problem with understanding lay, and
being able to supply the insight necessary to overcome the difficulty.'

\emph{Steve Noble has written}:
\lq In the 1990s the combinatorial group in Oxford was much smaller than now, comprising mainly Dominic and Colin McDiarmid and their graduate students. 
The Combinatorial Theory seminar was the main event of the week. 
We were fortunate to hear 
eminent speakers talking about topics 
relevant to our problems. 
But as we knew very little, much of what we listened to was largely incomprehensible. 
I was reassured when Dominic, with characteristic modesty, told me that he didn't understand much of a particularly lengthy and difficult proof. The 
highlight of the seminar was 
the discussion afterwards,
when Dominic would interrogate the speaker in a jovial but probing way. A speaker was never allowed to get away with a vague claim, and Dominic was never afraid to ask a question. Some of his questions were deep, and some less so, but there was an important third category that sounded silly at first but somehow pointed to the nub of a topic in a way which no one else had considered, often suggesting connections between apparently unrelated areas.
I assumed that 
it was always like this, but I have never come across anyone who could do whatever it was Dominic did to make the discussion so lively and informative. I learnt a huge amount from Dominic, much of it from these informal discussions.\rq

\emph{Colin McDiarmid has written}:
\lq Before I came to Oxford to start graduate work I had certainly enjoyed mathematics, but with Dominic as supervisor the game changed.  Mathematics became very exciting, with a playful side, something to share.  From the beginning Dominic and Bridget made me feel very much at home in Merton and in Oxford.  

In my first meeting with Dominic we walked in Christ Church Meadows in sight of Merton; he often told the story that my first question was not about mathematics or `where is the Merton library?’ but `where are the Merton rugby pitches?’.  He was a great story teller, always with a smile and a laugh.  Another early meeting was at Blenheim, with fellow new graduate Frank Dunstan 
$\dots$ and a  frisbee.  Often we would meet at his home, with family around.   
More happy memories are at workshops or conferences where it was always so full of life around Dominic.  
He was a great supervisor, certainly my model.  Before long he was a good friend as well as a mentor, then a much respected and loved colleague, and now much missed.\rq

\section{Work and influence}

It may be said that Dominic's mathematical trajectory was influenced most
by John Hammersley and Bill Tutte --- the former imparted  a
love of problems associated with counting and chance, 
and the work of the latter was foundational to Dominic's interests in combinatorics.
The human angle of mathematics as a participatory activity played a  key role for him.
He thrived off discussions with colleagues, usually (ex)-students, throwing ideas and often
wild questions around, possibly over tea in his home or in the Oxford Mathematical Institute.
One thing would lead to another, and papers, books, and dinners flowed, usually in the company
of members of the community that circulated around and were inspired by such meetings.

The 1960s and \rq 70s were a glorious period for combinatorial theory, with Dominic active at the heart of
the development in the UK. He organised and edited the proceedings of the 1969 Oxford Conference that came to 
be viewed as the First British Combinatorial  Conference (BCC); the list of 35 speakers included 
Paul Erd\H os, Richard Guy, Mark Kac, and Roger Penrose.
This gave birth to the continuing series of BCC meetings (\cite{NB}) and later the British Combinatorial
Committee of which Dominic was the Chair for the period 1983--1987.

The remainder of this memoir is a brief overview of Dominic's contributions to probability and combinatorics.  
Section \ref{sec:perc} is centred around his work on percolation and self-avoiding walk, 
largely done as a D.Phil.\ student;
the results summarised there have had great 
impact over the intervening 60 years. Sections \ref{sec:mat} and \ref{sec:comp} are directed towards his work,
as an individual and often jointly with his students, on matroids and complexity.  

Dominic began taking postgraduate students very soon after his own doctoral
graduation, beginning with Adrian Bondy
who graduated in 1969.
On his academic family tree \cite{DW2} are listed 28 students, with 235 scientific offspring 
in all at the time of this memoir. 
\emph{Graham Farr, Dillon Mayhew, and James Oxley have written in \cite{FMO}}:
\lq Dominic was a very effective supervisor of research students. 
[$\dots$] He was flexible in his approach and adept at finding a productive mix of patience, firmness, encouragement, plain speaking, and inspiration. [$\dots$] Time and again he brought out the best [$\dots$], inspiring enduring appreciation and affection.'

In his books, Dominic achieved a high level of communication (and excellent Amazon reviews)
through clear exposition and a minimum of mathematical prerequisites. In addition to his more
elementary volume \cite{MR3243603} on probability (which, wrote 
one eminent reviewer, reads as though it was written in a punt), 
he wrote more advanced texts on \emph{Matroid Theory} \cite{MR0427112},
\emph{Codes and Cryptography} \cite{MR0959137},  \emph{Complexity: Knots, Colourings and Counting} 
\cite{MR1245272},  and
\emph{Complexity and Cryptography} (with John Talbot) \cite{MR2221458}.

Dominic's friends and colleagues marked his retirement from Oxford in 2005 
by the publication of the collection \cite{ccc} of articles on topics close to his heart.
Moreover, he was awarded in 2006 the degree of Doctor of Mathematics, \emph{honoris causa}, at the University of 
Waterloo in recognition of his work in combinatorial mathematics.
A full list of Dominic Welsh's publications is available at \cite{djaw-bib}.

\section{Percolation and self-avoiding walks}\label{sec:perc}

Dominic's D.Phil.\ supervisor, John Hammersley, was a pioneering mathematician who contributed
some of the very best early work on important topics including
self-avoiding walks and percolation.
Each of
their two joint papers has had substantial impact on fields that have grown 
enormously in visibility over the decades.   

\subsection{First-passage percolation}

Percolation is the canonical model for a random spatial medium.  It is usually
viewed as a \emph{static} model in the sense that there is no dependence on time and no evolutionary aspect.
The \emph{first-passage percolation} process, introduced by Hammersley and Welsh in \cite{MR0198576} in 1965,
includes a time variable, and thus leads to beautiful questions involving growing random sets and their
asymptotic behaviour. This work formed the core of Dominic's D.Phil\ thesis \cite{DPhil-thesis}, and
has had great influence over the intervening 60 years on the theory of random spatial growth.

Amongst the major advances of this work is the introduction of the concept of a \emph{subadditive} stochastic process,
that is, a stationary family $(X_{s,t}: 0\le s<t <\oo )$ of random variables satisfying
subadditivity: for $s<u<t$, we have $X_{s,t}\le X_{s,u}+X_{u,t}$. The authors proved an early version
of the now famous subadditive ergodic theorem, which states that, subject to a suitable moment condition, the limit
$\lim_{t\to\oo} X_{0,t}/t$ exists. This theorem has been improved and optimised
in various ways since 1964 by Kingman \cite{jk1,jk2}, Liggett \cite{tml}, and others;
Kesten and Hammersley weakened the subadditive condition to a distributional assumption, and so on.
The subadditive ergodic theorem is now one of the major tools of probability theory.

\subsection{Self-avoiding walks}

A \emph{self-avoiding walk} (SAW) on a graph $G$ is a path that visits no vertex twice or more. The
basic SAW problem is to estimate the number $\sigma_n(G)$ of distinct $n$-step SAWs starting at a given vertex.
This problem is fundamental to the theory of polymers as studied by Flory
and co-workers, \cite{flory}. For concreteness (and a reason that will soon be clear), 
let $G$ be the hexagonal lattice, denoted $\HH$. It is now standard that $\sigma_n$ is a 
(deterministic) subadditive sequence
in that $\sigma_{m+n} \le \sigma_m + \sigma_n$, and it follows immediately that the limit
$\kappa=\lim_{n\to\oo} \sigma_n^{1/n}$
exists.
The constant $\kappa=\kappa(\HH)$ is called the \emph{connective constant} of $\HH$, and one has by subadditivity that
$\sigma_n \ge \kappa^n$. In \cite{MR0139535}, Hammersley and Welsh sought upper bounds for $\sigma_n$, and were
able to prove that there exists $\gamma$ such that $\sigma_n \le \kappa^n\gamma^{\sqrt n}$.
This they achieved by showing that a SAW can be decomposed as a union of so-called bridges, 
and by observing that counts of bridges are superadditive rather than subadditive.
(Actually they worked with the $d$-dimensional cubic lattice, but their argument is largely valid 
in greater generality.)

Two observations are made about the impact of this work. Firstly, it is believed (but not yet proved) that
the true correction term in the two-dimensional case is a power of $n$ rather than an exponential of $\sqrt n$. 
More precisely, it is believed, for a $d$-dimensional lattice $\sL$, that there exist constants $A=A_d$ and 
$\gamma=\gamma_d$ such that
\begin{equation}\label{eq:saw}
\sigma_n \sim A n^{\gamma-1}\kappa(\sL)^n,
\end{equation}
(subject to a logarithmic correction when $d=4$)
and moreover $\gamma_2=\frac{43}{32}$. 
Equation \eqref{eq:saw} was proved in 1992 by Hara and Slade for 
the cubic lattice in $5$ and more dimensions, in which case
we have $\gamma=1$ (see \cite{MR2986656}).
The case of two dimensions seems especially hard, and the best rigorous work so far
appears to be that of \cite{MR4124522} where the $\gamma^{\sqrt n}$ is replaced by 
$\gamma^{n^{\frac12-\eps}}$ for some $\eps>0$.    
The full picture in two dimensions may emerge only when we have a proper understanding of
the relationship between SAWs and the stochastic Loewner evolution process SLE$_{8/3}$.
 
Returning to the case of the hexagonal lattice $\HH$,
the biggest result on SAWs in recent years is the exact calculation of $\kappa(\HH)$ 
by Duminil-Copin and Smirnov \cite{MR2912714}, namely
$\kappa(\HH) = \sqrt{2+\sqrt 2}$,
as conjectured by Nienhuis using conformal field theory. The proof is a combination
of a deeply original argument combined  (in its \lq easier' part) with the bridge arguments of
Hammersley and Welsh \cite{MR0139535}.

\subsection{Russo--Seymour--Welsh theory for percolation}

Dominic loved problems, and one of his favourites in the 1970s
was to prove that the
critical probability of bond percolation on the square lattice $\ZZ^2$
equals $\frac12$. This captivating and beautiful conjecture
seemed to defy all attempts. Dominic and Paul Seymour constructed a new tool now named
the RSW lemma after Russo \cite{MR0488383} and Seymour--Welsh \cite{MR0494572}.

Bond percolation on $\ZZ^2$ is given as follows. Let $0<p<1$, and declare each edge of $\ZZ^2$
\emph{open} with probability $p$ and \emph{closed} otherwise (with independence between edges).
Let $\theta(p)$ be the probability that the origin lies in an infinite open component, and define the 
\emph{critical probability}
$\pc=\sup\{p: \theta(p)=0\}$.
Hammersley conjectured that $\pc=\frac12$, and Ted Harris \cite{MR0115221} proved that $\pc\ge \frac12$.
The conjecture is supported by a symmetry between open edges of $\ZZ^2$ and closed edges of its dual graph. 

Dominic's idea was to study what he called \lq sponge percolation', by which he meant the probability $\pi_p({m,n})$
that an $m\times n$ rectangle is crossed from left to right by an open path. 
Together with Seymour, he proved by a 
complicated geometrical and probabilistic construction that there exists a continuous function $f:(0,1)\to (0,1)$ such that
$\pi_p(4m,2m) \ge f(\pi_p(2m,2m))$ for all $m$.
This can be read as saying that, if there is a decent probability of a crossing of a square,
 then there is a decent probability of a crossing of a rectangle. This in turn implies a lower bound for 
 the probabilities of still longer paths and cycles.  
 
 This RSW lemma is one of the fundamental tools of percolation and related topics (see \cite{MR4557761}
 for a recent result of RSW type). It is the basic tool used by Harry Kesten in his celebrated
 proof that indeed $\pc(\ZZ^2)=\frac12$, \cite{MR0575895}.

\section{Matroid theory}\label{sec:mat}

A systematic study of the concept of linear independence was pursued in the 1930s through the work of van den Waerden,
Whitney, and others, and it was Whitney who coined the word \lq matroid' in his fundamental paper \cite{whitney}.
Matroids attracted the attention
of Tutte, Rado, and others in the 1950s, and their study grew gradually within the combinatorial community.
They have gained  greater prominence in recent years, with connections to a number of fields
beyond combinatorics including topology, geometric model theory, non-commutative geometry, and more recently
algebraic geometry.  

Dominic's interests in combinatorics were broad, with matroids as a principal theme.
He was brought to matroids by a \lq\lq stimulating seminar on the applications of matroids\rq\rq\
by Crispin Nash-Williams in 1966 (see \cite[p.\ v]{MR0427112}).
He played a major role in the
systemisation of matroid theory and its associations over the next decade.
He enriched the field with numerous theorems and conjectures, and the latter have often proved 
fruitful.
He wrote 31 papers between 1967 and 2013 with the word \lq matroid' in the title, in addition to
his monograph \emph{Matroid Theory} of 1976.

In some of his earlier papers, Dominic pursued results inspired by graph theory; he extended
to matroids 
Kruskal's greedy algorithm for constructing a minimal spanning tree in a weighted graph
(this result is ascribed in \cite{ox-ccc} to a pioneering  paper of Bor\r{u}vka from 1926
written in Czech, see \cite{MR1825599});
he proved a version of Euler's theorem for binary matroids. He worked also on transversal theory;
he investigated  the relationship between Hall's marriage theorem and Rado's condition for the existence of a transversal;
with Joan de Sousa he proved that every binary transversal matroid is graphic.
He advanced a number of other problems, including  the representability of matroids, and their relationship to block designs.

Enumeration provided another theme of study, beginning in 1969 with a constructive proof of Crapo's lower
bound on the number $g(n)$ of matroids on a set of $n$ elements. 
With Michael Piff, he followed this by improving the lower bound to $\log_2 g(n) \ge 2^n n^{-\frac52} \log_2n$
for large $n$.
Knuth later \cite{MR0335312} replaced $\frac52$ by $\frac32$ by showing that
$\log_2 (n!\, g(n)) \ge \binom n {n/2}/(2n)$, and it turns out that this bound is, in a sense, not too far from the 
correct asymptotics for $g(n)$ (see \cite{bansal}).

A significant number of counting problems rotate around the related concepts of unimodality and log-concavity.
A real sequence $(w_i: 1\le i\le n)$ is \emph{unimodal} if there exists $k$ such that 
$w_1\le w_2 \le\dots\le w_{k-1}\le  w_k\ge w_{k+1}
\ge\dots\ge w_n$, and \emph{log-concave} if $w_i^2\ge w_{i-1}w_{i+1}$ for all $i$. 
Dominic asked whether (i) the number $f_i(n)$ of rank-$i$ non-isomorphic 
matroids on an $n$-set constitutes a unimodal sequence, and (ii) the number $i_k$ of $k$-element independent sets
of a matroid is unimodal. A conjecture of unimodality may be extended to a stronger (for positive terms)
conjecture of log-concavity
(\cite{MR0349445}).

A further question is concerned with the coefficients of the characteristic polynomial of a matroid $M$ on the set $E$,
$$
\chi_M(\lambda)= \sum_{I \subseteq E} (-1)^{|I|}\lambda^{r(E)-r(I)},
$$
where $r$ denotes rank; this may be viewed as the matroid version of the chromatic 
polynomial of a graph.
 Excepting trivial cases, the coefficients of this power series oscillate in sign.
True to form, in \cite[Exerc.\ 15.3.5]{MR0427112}
Dominic invited the reader to prove the log-concavity of the sequence  of these coefficients
(earlier related conjectures are due to Gian-Carlo Rota and Dominic's student Andrew Heron); 
this exercise is marked delicately with a small \lq o' to indicate an open problem.
Forty years on, this, and the above log-concavity of the counts of independent sets,
were proved by Adaprasito, Huh, and Katz \cite{MR3862944} in the third  
of a series of related works by June Huh and his co-authors.
Their proof is based on a novel Hodge theory for matroids. The reader is referred to \cite[Sec.\ 15.2]{MR2849819} for 
an account of three log-concavity conjectures for matroids, and to \cite{MR4680241} for an overview
of the work of \cite{MR3862944} and its precursors.  

Fuller accounts of Dominic's work on matroids are found in \cite{FMO} and \cite{ox-ccc}.

\section{Complexity theory}\label{sec:comp}

\subsection{Computational complexity}

Students of Hammersley tended to have computational tendencies, and Dominic was no exception. 
His time at Bell Laboratories was influential, and he soon developed his own ideas for computational
methods.  His work as a D.Phil.\ student inspired two letters about PERT networks
to the editor of {\em Operations Research}, and  
he followed this with Martin Powell in proposing the so-called Welsh--Powell algorithm \cite{wp} for colouring a graph.
This is a sequential colouring algorithm: first list the vertices in decreasing order of degree; then construct 
greedily an independent set containing the first vertex, and colour this set with the first colour; iterate. 

Dominic's interests in complexity developed in the 1980s, inspired in part by the volume 
\cite{MR0519066} of Garey and Johnson. He wrote a number of influential review articles,
and encouraged his students to work on  the classification of
graph-theoretic complexity problems such as thickness and colouring. 
The following extract from the MathSciNet review
of his survey \cite{MR0671913}
with his student Tony Mansfield might be written equally about any of his expository work: 
`The paper is well written and is eminently readable. [$\dots$] 
The authors have artfully employed the best of both intuition and formalism to achieve succinctness and clarity.\rq

\subsection{Tutte polynomials}

Bill Tutte introduced his eponymous polynomial in 1947 (see \cite{handbook34});
its prominence has since grown vastly (see \cite{handbook}).

The Tutte polynomial of a finite graph $G=(V,E)$ is given by
$$
T_G(x,y)= \sum_{A\subseteq E} (x-1)^{r(E)-r(A)}(y-1)^{|A|-r(A)},
$$
where $r$ denotes rank. The study of such polynomials 
is now a major part of graph (and matroid) theory. By choosing $(x,y)\in\CC^2$ suitably, one obtains counts of
 important features of $G$ such as spanning trees and forests. 
In addition, $T_G$ includes the chromatic, flow, and knot polynomials as well
as the Whitney rank generating function and the random-cluster/Potts partition functions. 

\begin{SCfigure}
\centering
\caption{The real parameter plane $\RR^2$ of the Tutte polynomial $T_G$. The red hyperbola $H_1$
is given by $(x-1)(y-1)=1$. The blue hyperbola $H_2$, $(x-1)(y-1)=2$, corresponds to the partition function 
of the Ising model on $G$ (when planar). This was probably Dominic's favourite figure during the 1990s.}
\q\includegraphics[width=0.5\hsize]{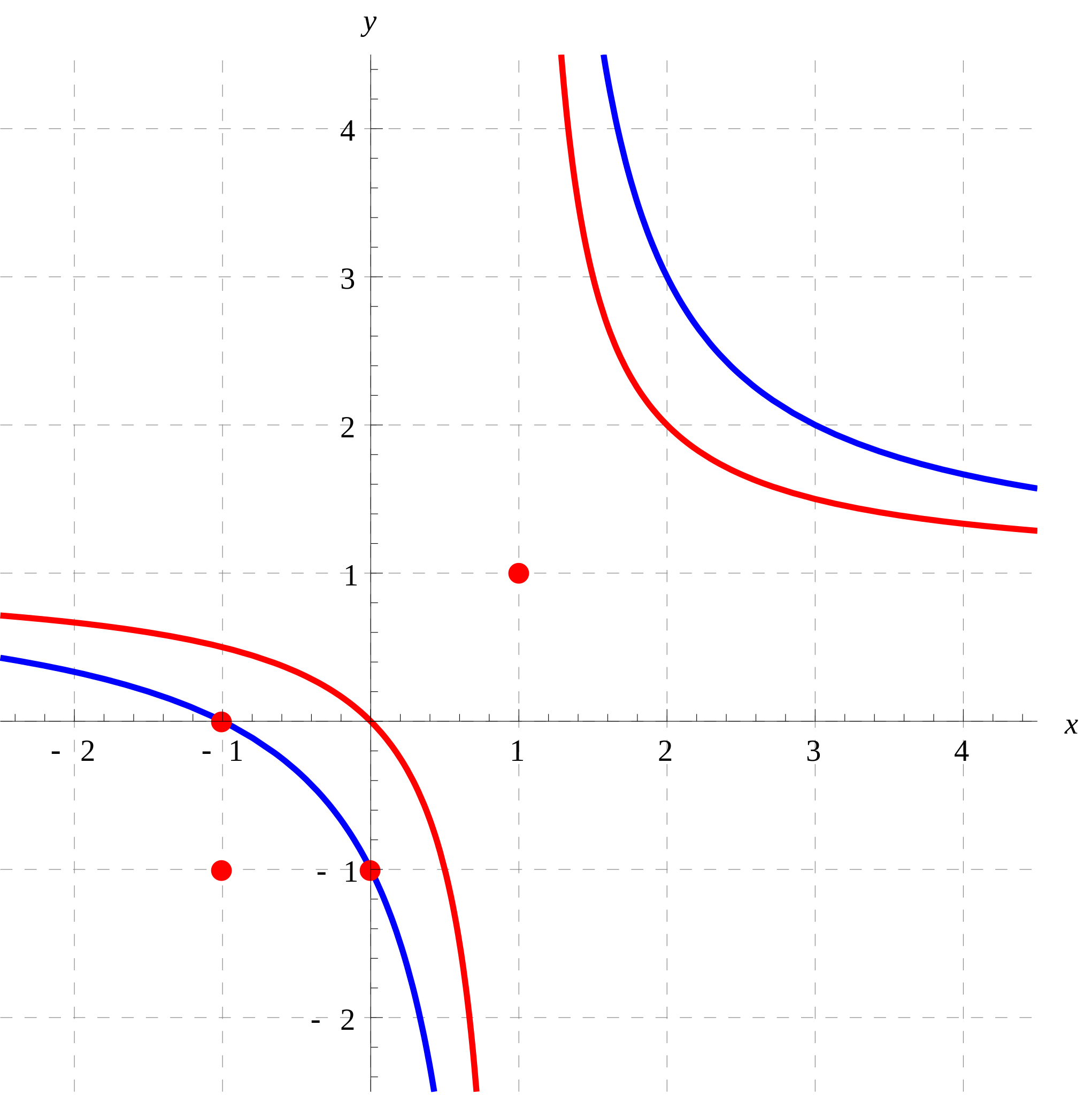}
\label{fig:comp}
\end{SCfigure}

Dominic devoted much of his later years to studying Tutte polynomials. In an important
article \cite{MR1049758}, written with Jaeger and Vertigan,
he initiated a systematic study of the complexity of calculating $T_G(x,y)$ for certain classes of
graph $G$ and for different values of $(x,y)\in\CC^2$ (his book \cite{MR1245272}
has been influential in this field). Let $R$ be the set of red points of Figure \ref{fig:comp}
(for simplicity, we consider only the real case $(x,y)\in\RR^2$).
They proved that $T_G(x,y)$ is computable in polynomial time for $(x,y)\in R$, 
whereas its computation is $\#$P-hard on $\RR^2\setminus R$.

Moreover, for \emph{planar} graphs, $T_G$ may be computed in polynomial time on the hyperbola $H_2$ of
Figure \ref{fig:comp}; this follows 
by the classic representation by Kasteleyn and Fisher of the Ising partition function in terms of dimer
configurations (also known as complete matchings).

This work was the first of a series of papers by Dominic's students and others on the complexity of
the Tutte polynomial and its cousins. In addition, it led to extensions
of the work to matroids, and also to the problem of approximating Tutte polynomials. 

A \emph{fully polynomial randomised approximation scheme} (FPRAS)  
for a function $f(x)$ is a randomised algorithm
that, for all $\eps>0$ and all instances $x$, 
outputs a random value $\what f(x,\eps)$ satisfying 
$$
\PP\bigl(\what f(x,\eps)/f(x) \in(1-\eps,1+\eps)\bigr)\ge \tfrac34,
$$
and that runs in time polynomial in $|x|$ and $1/\epsilon$.
For what graphs $G$ and $(x,y)\in\RR^2$ does there exist an FPRAS to compute  $T_G(x,y)$?
Dominic's student James Annan considered the case of dense graphs: for $\alpha>0$,  $\sG_\alpha$ is the class of 
$\alpha$-\emph{dense} graphs $G=(V,E)$ satisfying $|E| \ge \alpha|V|$. 
With Noga Alon and Alan Frieze, Dominic showed in \cite{MR1368847} the existence of an FPRAS when $x,y\ge 1$
for the class $\sG_\alpha$ with $\alpha>0$ (with a stronger conclusion when $\alpha>\frac12$). 
Never lacking in bravery, Dominic has conjectured that the condition on $\alpha$ may be removed, 
but the jury is still out on that (see \cite{handbook10}).

One further conjecture of Dominic deserves mention. With his student Criel Merino he conjectured in \cite{MR1772357}
that the number of spanning trees of a graph is no greater than the maximum of the number of acyclic orientations 
and the number of
totally cyclic orientations. This amounts to the inequality
$$
T_G(1,1)\le \max\{T_G(2,0), T_G(0,2)\},
$$
and this has been open since their 1999 paper.

In studying randomised algorithms, one encounters the problem of ascertaining the rate of convergence
of associated Markov chains. Dominic was thus led to the problem of finding
a way to simulate a random planar graph, and of identifying its properties.  
He initiated this project with Alain Denise
and Marcio Vasconcellos in \cite{MR1393702}, and developed it further with Colin McDiarmid and Angelika Steger in 
\cite{MR2117936}.
The required techniques turn out to be largely combinatorial and analytic, and it is fitting
that one of Dominic's projects in later life should be so close to the celebrated work \cite{tutte63} of Bill Tutte on map
enumeration from the time of the former's D.Phil.\ study.

\section*{Acknowledgements}

The author thanks friends and colleagues who have kindly contributed to and commented
on this memoir. He is grateful above all to Bridget Welsh for  
her kindnesses over many years, and for the wonderful conversations that have informed this account.
Graham Farr, Dillon Mayhew, and James Oxley have kindly contributed through their very useful account \cite{FMO} of Dominic's life and work; thanks also to Graham
for sharing the copy of Dominic's D.Phil.\ thesis supplied by the Bodleian Library.
The portrait of Dominic was taken by Bridget Welsh in 2006.
Figure \ref{fig:comp} was published on Wikipedia under the CC-BY-SA 3.0 licence.
Han Kimmett of Merton College kindly advised on the details of Dominic's association with Merton.

\providecommand{\bysame}{\leavevmode\hbox to3em{\hrulefill}\thinspace}
\providecommand{\MR}{\relax\ifhmode\unskip\space\fi MR }
\providecommand{\MRhref}[2]{%
  \href{http://www.ams.org/mathscinet-getitem?mr=#1}{#2}
}
\providecommand{\href}[2]{#2}

\newpage
\renewcommand\refname{Publications of Dominic Welsh}
\small

\newcommand\cl[1]{\smallskip\goodbreak\centerline{#1}\nobreak\smallskip\nobreak}

\end{document}